%% file: AB2008_v5.tex
\documentclass{llncs}
\usepackage{makeidx}
\usepackage{rotating}
  \usepackage{graphicx}
  \usepackage{psfrag}
  \usepackage{amsfonts}
\usepackage{amsmath,amssymb,color}
\usepackage{url}

  \usepackage{amsfonts}
\graphicspath{{figs/}}

\newcommand{\R}{\mathbb R}

\newcommand{\bigst}{\mathbin \big |}

\newcommand{\st}{\mathbin |}

\def\comment#1{\textit{[#1]}}
\def\comment#1{}

\begin{document}

\title{The Geometry of the Neighbor-Joining Algorithm for Small Trees}

\author{Kord Eickmeyer\inst{1} and Ruriko Yoshida\inst{2}}

\institute{Institut f\"ur Informatik, Humboldt-Universit\"at zu Berlin, Berlin, Germany \and University of Kentucky, Lexington, KY, USA}

\maketitle

\begin{abstract}
In 2007, Eickmeyer et al. showed that 
the tree topologies outputted by the Neighbor-Joining (NJ) algorithm and the 
balanced minimum evolution (BME) method for phylogenetic reconstruction 
are each determined by a polyhedral subdivision of the space of dissimilarity 
maps ${\R}^{n \choose 2}$, where $n$ is the number of taxa.
In this paper, we will analyze the behavior of the Neighbor-Joining
algorithm on five and six taxa and study the geometry and 
combinatorics of the polyhedral subdivision of the space of dissimilarity maps
for six taxa as well as hyperplane representations of each polyhedral 
subdivision.
We also study simulations for one of the questions stated by 
Eickmeyer et al., that is, the robustness of the NJ algorithm to 
small perturbations of tree metrics, with tree models which are known to
be hard to be reconstructed via the NJ algorithm.
\end{abstract}
\section{Introduction}

\begin{figure}[htp]
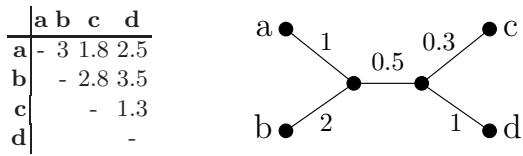

\vskip -0.2in
\begin{center}
  \parbox[c]{3.2cm}{
    \begin{tabular}{r|cccc}
      & {\bf a} & {\bf b} & {\bf c} & {\bf d}
      \\
      \hline
	  {\bf a} & - & $3$ & $1.8$ & $2.5$ \\
	  {\bf b} &   &  -  & $2.8$ & $3.5$ \\
	  {\bf c} &   &     &   -   & $1.3$ \\
	  {\bf d} &   &     &       &   -
    \end{tabular}
  }
  \parbox[c]{3.5cm}{
    \resizebox{\linewidth}{!}{\input figs/sampletree.pstex_t}
  }
\end{center}
\caption{ The NJ algorithm takes a matrix of pairwise distances (left)
as input and computes a binary tree (right). If there is a tree such
that the distance matrix can be obtained by taking the length of the
unique path between two nodes, NJ outputs that tree.}\label{fig:sampletree}
\vskip -0.2in
\end{figure}

The Neighbor-Joining (NJ) algorithm was introduced by Saitou and Nei
\cite{Saitou1987} and is widely used to reconstruct a phylogenetic
tree from an alignment of DNA sequences because of its accuracy and
computational speed.  The NJ algorithm is a distance-based method
which takes all pairwise distances computed from the data as its
input, and outputs a tree topology which realizes these pairwise
distances, if there is such a topology (see Fig. \ref{fig:sampletree}). 
Note that the NJ algorithm is
consistent, i.e., it returns the additive tree if the input distance
matrix is a tree metric. \comment{ However, we usually estimate all pairwise
distances via the maximum likelihood estimation (MLE).} If the input
distance matrix is not a tree metric, then the NJ algorithm returns a
tree topology which induces a tree metric that is ``close'' to the
input.  Since it is one of the most popular methods for reconstructing
a tree among biologists, it is important to study how the NJ algorithm
works.

A number of attempts have been made to understand the good results
obtained with the NJ algorithm, especially given the
problems with the inference procedures used for estimating pairwise
distances.  For example, Bryant showed that the {\em
Q-criterion} (defined in \eqref{eqn:qcrit} in Section
\ref{qcriterion}) is in fact the unique selection criterion which
is linear, permutation invariant, and {\em consistent}, i.e., it
correctly finds the tree corresponding to a tree metric \cite{Bryant2005}.
Gascuel and Steel gave a nice review of how the NJ algorithm works \cite{Steel2006}.

One of the most important questions in studying the behavior of the NJ
algorithm is to analyze its performance with pairwise distances that
are not tree metrics, especially when all pairwise distances are
estimated via the maximum likelihood estimation (MLE).  
In 1999, Atteson showed that
if the distance estimates are at most half the minimal
edge length of the tree away from their true value then the NJ algorithm will reconstruct the correct tree \cite{Atteson99}.  However, Levy et al.
noted that Atteson's criterion frequently fails to be
satisfied even though the NJ algorithm returns the true tree topology \cite{Levy2005}.
Recent work of \cite{Mihaescu2006} extended Atteson's work.
Mihaescu et al. showed that the NJ algorithm returns the true tree
topology when it works locally for the quartets in the tree \cite{Mihaescu2006}.  This
result gives another criterion for when then NJ algorithm returns the
correct tree topology and Atteson's theorem is a special case of
Mihaescu et al.'s theorem. \comment{ The results of \cite{Atteson99} and
\cite{Mihaescu2006}, however, depend on a particular input data and to analyze
the performance of the NJ algorithm for {\em any} input data we need
to study the sample space ($\R^{{n \choose 2}}$, where $n$ is
the number of taxa).  Therefore, this justifies the study of sample space which
is a subspace of the vector space $\R^{{n \choose 2}}$, where $n$ is
the number of taxa.  }

For every input distance matrix, the NJ algorithm returns a certain
tree topology. It may happen that the minimum Q-criterion is taken by
more than one pair of taxa at some step. In practice, the NJ algorithm
will then have to choose one cherry in order to return a definite
result, but for our analysis we assume that in those cases the NJ
algorithm will return a set containing all tree topologies resulting
from picking optimal cherries.
There are only finitely many tree topologies, and for every topology
$t$ we get a subset $D_t$ of the sample space (input space) such that for all
distance matrices in $D_t$ one possible answer of the NJ algorithm is
$t$. We aim at describing these sets $D_t$ and the relation between
them.
One notices that the Q-criteria are all linear in the pairwise distances.  
The NJ algorithm will pick 
cherries in a particular order and output a particular tree $t$ if and only 
if the pairwise distances satisfy a system of linear inequalities, whose 
solution set forms a polyhedral cone in  ${\R}^{n \choose 2}$.  
Eickmeyer et al. called such a cone a 
{\em Neighbor-Joining cone,} or {\em NJ cone}.  Thus the NJ algorithm will
output a particular tree $t$ if and only if the distance data lies in a 
union of NJ cones \cite{NJME}.

In \cite{NJME}, Eickmeyer et al. studied the optimality of the NJ 
algorithm compared to
the balanced minimum evolution (BME) method and focused on polyhedral 
subdivisions of the space of dissimilarity maps for the BME method and 
the NJ algorithm.  Eickmeyer et al. also studied the 
geometry and combinatorics of the NJ cones for $n = 5$ in addition to the
{\em BME cones} for $ n \leq 8$.
Using geometry of the NJ cones for $n = 5$, they showed 
that the polyhedral 
subdivision of the space of dissimilarity maps with the NJ algorithm
does not form a fan for $n \geq 5$ and that the union of the 
NJ cones for a particular tree topology is not convex.  
This means that the NJ algorithm is not convex, i.e.,  there are distance 
matrices $D, \, D'$, such that NJ produces the same tree $t_1$  
on both inputs $D$ and $D'$, but it produces a 
different tree 
$t_2 \neq t_1$ on the input $(D + D')/2$ 
(see \cite{NJME} for an example).

In this paper, we focus on 
describing geometry and combinatorics of the NJ cones for six taxa
as well as some simulation study using the NJ cones for one of the questions in 
\cite{NJME}, that is, what is the robustness of the NJ algorithm to 
small perturbations of tree metrics for $n = 5$.
This paper is organized as follows: In Section \ref{sec:nj} we will describe the NJ algorithm and define the NJ cones.  Section \ref{C5} states the hyperplane
representations of the NJ cones for $n = 5$.  Section \ref{6taxa} describes 
in summary the geometry and combinatorics of the NJ cones for $n = 6$.
Section \ref{simulations} shows some simulation studies on the robustness of 
the NJ algorithm to small perturbations of tree metrics for $n = 5$ with
the tree models from \cite{Ota2000}.  We end by discussing some open problems
in Section \ref{problem}.
\vskip -0.2in

\section{The Neighbor-Joining algorithm}
\label{sec:nj}

\subsection{Input data}
The NJ algorithm is a {\em distance-based method} which takes a {\em
distance matrix}, a symmetric matrix $(d_{ij})_{0\leq i, j\leq n-1}$
with $d_{ii} = 0$ representing pairwise distances of a set of $n$ taxa
$\{0, 1, \ldots, n - 1\}$, as the input.  Through this paper, we do
not assume anything on an input data except it is symmetric and $d_{ii}
= 0$. Because of symmetry, the input can be seen as a vector of
dimension $m := \binom{n}{2} = \frac{1}{2}n(n-1)$. We arrange the
entries row-wise. 
We denote row/column-indices by pairs of letters such as $a$, $b$, $c$, $d$, while denoting single indices into the ``flattened'' vector by letters $i, j, \dots$. The two indexing methods are used simultaneously in the hope that no confusion will arise. Thus, in the four taxa example we have $d_{0,1} = d_{1,0} = d_0$. In general, we get $d_i = d_{a,b} = d_{b,a}$ with
\begin{equation*}
a = \max\left\{k \bigst \frac{1}{2}k(k-1) \leq i\right\} = 
\left\lfloor \frac{1}{2}+ \sqrt{\frac{1}{4}+2i} \right\rfloor,
b = i-\frac{1}{2}(a-1)a,
\end{equation*}
and for $c>d$ we get
\begin{equation*}
d_{c,d} = d_{c(c-1)/2 + d}.
\end{equation*}
\vskip -0.1in

\subsection{The Q-Criterion}\label{qcriterion}

The NJ algorithm starts by computing the so called \emph{Q-criterion} or the {\em cherry picking criterion}, given by the formula
\begin{equation}
\label{eqn:qcrit}\tag{Q}
q_{a,b} := (n-2)d_{a,b} - \sum_{k=0}^{n-1} d_{a,k} - \sum_{k=0}^{n-1} d_{k,b}.
\end{equation}
This is a key of the NJ algorithm to choose which pair of taxa is a neighbor.
\begin{theorem}[Saitou and Nei, 1987, Studier and Keppler, 1988 \cite{Saitou1987,Studier1988}]
\label{Saitou}
Let $d_{a,b}$ for all pair of taxa $\{a, b\}$ be the tree metric corresponding to the tree $T$.  Then the pair $\{x, y\}$ which minimizes $q_{a,b}$ for all pair of taxa $\{a, b\}$ forms a neighbor.
\end{theorem}

Arranging the Q-criteria for all pairs in a matrix yields again a
symmetric matrix, and ignoring the diagonal entries we can see it as a
vector of dimension $m$ just like the input data. Moreover, the
Q-criterion is obtained from the input data by a linear
transformation:
\begin{equation*}
\mathbf{q} = A^{(n)}\mathbf{d},
\end{equation*}
and the entries of the matrix $A^{(n)}$ are given by
\begin{equation}
\label{eqn:adef}
A^{(n)}_{ij} = A^{(n)}_{ab,cd} = \begin{cases}
n-4 & \text{if }i = j,\\
-1  & \text{if }i \not = j \text{ and }\{a,b\}\cap \{c,d\} \not = \emptyset,\\
0   & \text{else},
\end{cases}
\end{equation}
where $a > b$ is the row/column-index equivalent to $i$ and likewise for $c > d$ and $j$. When no confusion arises about the number of taxa, we abbreviate $A^{(n)}$ to $A$.  

After computing the Q-criterion $\mathbf{q}$, the NJ algorithm
proceeds by finding the minimum entry of $\mathbf{q}$, or, equivalently, the
maximum entry of $-\mathbf{q}$. The two nodes forming the chosen pair
(there may be several pairs with minimal Q-criterion) are then joined
(``cherry picking''), i.e., they are removed from the set of nodes and
a new node is created.  Suppose out of our $n$ taxa
$\{0,\ldots,n-1\}$, the first cherry to be picked is $m-1$, so the
taxa $n-2$ and $n-1$ are joined to form a new node, which we view as
the new node number $n-2$. The reduced pairwise distance matrix is one
row and one column shorter than the original one, and by our choice of
which cherry we picked, only the entries in the rightmost column and
bottom row differ from the original ones. Explicitly,
\begin{equation*}
{d'}_i =
\begin{cases}
d_i & \text{for }0 \leq i < \binom{n-2}{2}
\\
\frac{1}{2}(d_i + d_{i+(n-2)} - d_{m-1})& \text{for }\binom{n-2}{2} \leq i < \binom{n-1}{2} 
\end{cases}
\end{equation*}
and we see that the reduced distance matrix depends linearly on the original one:
\begin{equation*}
\mathbf{d'} = R\mathbf{d},
\end{equation*}
with $R = (r_{ij}) \in \R^{(m-n+1)\times m}$, where
\begin{equation*}
r_{ij} = \begin{cases}
1   & \text{for }0 \leq i=j < \binom{n-2}{2}\\
1/2 & \text{for }\binom{n-2}{2} \leq i < \binom{n-1}{2}, j = i\\
1/2 & \text{for }\binom{n-2}{2} \leq i < \binom{n-1}{2}, j = i+n-2\\
-1/2 & \text{for }\binom{n-2}{2} \leq i < \binom{n-1}{2}, j = m-1\\
0 & \text{else}
\end{cases}
\end{equation*}
The process of picking cherries is repeated until there are only three
taxa left, which are then joined to a single new node.

We note that since new distances $d'$ are always linear combinations of 
the previous
distances, all Q-criteria computed throughout the NJ algorithm are linear 
combinations of the original pairwise distances.  Thus, for every possible 
tree topology $t$ outputted by the NJ algorithm (and every possible ordering $\sigma$ of picked cherries that results in topology $t$), there is a polyhedral cone $C_{T,\sigma} \subset \R^{n \choose 2}$ of dissimilarity maps.  The NJ algorithm will output $t$ and pick cherries in the order $\sigma$ iff 
the input lies in the cone $C_{T,\sigma}$.  We call the cones $C_{T,\sigma}$ {\em Neighbor-Joining cones}, or {\em NJ cones}.
\vskip -0.1in

\subsection{The shifting lemma}
\label{sec:shifting}

We first note that there is an $n$-dimensional linear subspace of
$\R^m$ which does not affect the outcome of the NJ algorithm (see
\cite{Mihaescu2006}). For a node $a$ we define its \emph{shift
vector} $\mathbf{s}_a$ by
\begin{equation*}
(\mathbf{s}_a)_{b,c} := \begin{cases}
1 & \text{if }a \in \{b,c\}
\\
0 & \text{else}
\end{cases}
\end{equation*}
which represents a tree where the leaf $a$ has distance 1 from all
other leaves and all other distances are zero. The Q-criterion of any
such vector is $-2$ for all pairs, so adding any linear combination of
shift vectors to an input vector does not change the relative values
of the Q-criteria. Also, regardless of which pair of nodes we join,
the reduced distance matrix of a shift vector is again a shift vector
(of lower dimension), whose Q-criterion will also be constant. Thus,
for any input vector $\mathbf{d}$, the behavior of the NJ algorithm on
$\mathbf{d}$ will be the same as on $\mathbf{d} + \mathbf{s}$ if
$\mathbf{s}$ is any linear combination of shift vectors. We call the
subspace generated by shift vectors $S$.

We note that the difference of any two shift vectors is in the kernel
of $A$, and the sum of all shift vectors is the constant vector with
all entries equal to $n$. If we fix a node $a$ then the set
\begin{equation*}
\{\mathbf{s}_a - \mathbf{s}_b \st b \not = a \}
\end{equation*}
is linearly independent. 
\vskip -0.1in

\subsection{The first step in cherry picking}
\label{sec:firststep}
After computing the Q-criterion $\mathbf{q}$, the NJ algorithm proceeds by finding the minimum entry of it, or, equivalently, the maximum entry of $-\mathbf{q}$. The set $cq_i \subset \R^m$ of all q-vectors for which $q_i$ is minimal is given by the normal cone at the vertex $-e_i$ to the (negative) simplex
\begin{equation*}
\Delta_{m-1} = \mathrm{conv}\{ -e_i \st 0\leq i \leq m-1 \}\subset \R^{m},
\end{equation*}
where $e_0,\ldots,e_{m-1}$ are the unit vectors in $\R^{m}$. The normal cone is defined in the usual way by
\begin{equation}
\label{eqn:cone}
\begin{split}
\mathcal{N}_{\Delta_{m-1}}(-e_i) :=&
\left\{
\mathbf{x} \in \R^m \st (-e_i,\mathbf{x}) \geq (\mathbf{p},\mathbf{x})\text{ for }\mathbf{p}\in\Delta_{m-1}
\right\}
\\
=&
\left\{
\mathbf{x} \in \R^m \st (-e_i,\mathbf{x}) \geq (-e_j,\mathbf{x})\text{ for }0\leq j \leq m-1
\right\},
\end{split}
\end{equation}
with $(\cdot,\cdot)$ denoting the inner product in $\R^m$.

Substituting $\mathbf{q} = A\mathbf{d}$ into \eqref{eqn:cone} gives
\begin{equation}
\label{eqn:conedef}
\begin{aligned}
\mathbf{q} \in cq_i
&\quad\Leftrightarrow\quad
i = \arg\max (-e_j,A\mathbf{d})
\\
&\quad\Leftrightarrow\quad
i = \arg\max (-A^T e_j,\mathbf{d})
\\
&\quad\Leftrightarrow\quad
i = \arg\max (-A e_j,\mathbf{d})
&&\text{because }A\text{ is symmetric}.
\end{aligned}
\end{equation}
Therefore the set $cd_i$ of all \emph{parameter} vectors $\mathbf{d}$ for which the NJ algorithm will select cherry $i$ in the first step is the normal cone at $-Ae_i$ to the polytope
\begin{equation}
\label{eqn:ppoly}
P_n := \mathrm{conv}\{-Ae_0,\ldots, -Ae_{m-1}\}.
\end{equation}
The shifting lemma implies that the affine dimension of the polytope $P_n$ is at most $m-n$. Computations using \texttt{polymake} show that this upper bound gives the true value. 

If equality holds for one of the inner products in this formula, then there are two cherries with the same Q-criterion.

As the number of taxa increases, the resulting polytope gets more complicated very quickly. By symmetry, the number of facets adjacent to a vertex is the same for every vertex, but this number grows following a strange pattern. See Table 1 for some calculated values.  We also computed f-vectors for $P_n$ via
{\tt polymake}.
With $n = 5$, we have $(1, 10, 45, 90, 75, 22, 1)$, with $n = 6$, 
$(1, 15, 105, 435, 1095, 1657, 1470, 735, 195, 25, 1)$, and
with $n \geq 7$ we ran {\tt polymake} over several hours and it took more than 
$9$GB RAM.  Therefore, we 
could not compute them.
\vskip -0.2in

\begin{table}
\label{tab:polytopes}
\begin{center}
\begin{tabular}{c|cccc}
\textbf{no. of} & \textbf{no. of}   & \textbf{dimension} & \textbf{facets through} & \textbf{no. of}\\
\textbf{taxa}   & \textbf{vertices} &                    & \textbf{vertex}         & \textbf{facets}\\
\hline
4  &   3 &  2  &      2 &      3\\
5  &  10 &  5  &     12 &     22\\
6  &  15 &  9  &     18 &     25\\
7  &  21 & 14  &    500 &    717\\
8  &  28 & 20  &    780 &  1,057\\
9  &  36 & 27  & 30,114 & 39,196\\
10 &  45 & 35  & 77,924 & 98,829\\
\end{tabular}
\end{center}
\caption{The polytopes $P_n$ for some small numbers of taxa $n$.}
\vskip -0.5in
\end{table}

\subsection{The cone $cd_i$}
Equation \eqref{eqn:conedef} allows us to write $cd_i$ as an intersection of half-spaces as follows:
\begin{equation}
\begin{split}
cd_i &= \left\{
\mathbf{x} \st (-A\mathbf{e}_i, \mathbf{x}) \geq (-A\mathbf{e}_j, \mathbf{x})
               \text{ for }j \not = i
\right\}
\\
&= \left\{
\mathbf{x} \st (-A(\mathbf{e}_i - \mathbf{e}_j), \mathbf{x}) \geq 0
               \text{ for }j \not = i
\right\}
\\
&=
\bigcap_{j\not = i}
\left\{
\mathbf{x} \st (-A(\mathbf{e}_i - \mathbf{e}_j), \mathbf{x}) \geq 0
\right\}
\end{split}
\end{equation}
We name the half-spaces, their interiors and the normal vectors defining them as follows:
\begin{equation}
\begin{split}
\mathbf{h}^{(n)}_{ij} &:= -A^{(n)}(\mathbf{e}_i - \mathbf{e}_j),
\\
H^{(n)}_{ij} &:= \left\{ \mathbf{x} \in \R^{m} \st ( \mathbf{h}^{(n)}_{ij}, \mathbf{x} ) \geq 0 \right\},
\\
{\mathring H}^{(n)}_{ij} &:= \left\{ \mathbf{x} \in \R^{m} \st ( \mathbf{h}^{(n)}_{ij}, \mathbf{x} ) > 0 \right\},
\end{split}
\end{equation}
where again we omit the superscript $(n)$ if the number of taxa is clear.

If there are more than four taxa, then this representation is not redundant: For any pair $i$ and $j$ of cherries, we can find a parameter vector $\mathbf{d}$ lying on the border of $H_{ij}$ (i.e., $(\mathbf{h}_{ij}, \mathbf{d}) = 0$) but in the interior ${\mathring H}_{ik}$ of the other half-spaces for $k \not = i, j$. One such $\mathbf{d}$ is given by
\begin{equation*}
d_k := \begin{cases}
2&\text{if }k = i\text{ or }k=j,\\
4&\text{else}.
\end{cases}
\end{equation*}

Thus we have found an $\mathcal{H}$-representation of the polyhedron $cd_i$ consisting of only $m-1$ inequalities. Note that Table 1 implies that a $\mathcal{V}$-representation of the same cone would be much more complicated, as the number of vectors spanning it is equal to the number of facets incident at the vertex $-A\mathbf{e}_i$.
\comment{
\begin{example}
In the case of four taxa, $P_3$ is a triangle with vertices
\begin{equation*}
p_0 = (0, 1, 1, 1, 1, 0)^\mathrm{T}, p_1 = (1, 0, 1, 1, 0, 1)^\mathrm{T},p_2 = (1, 1, 0, 0, 1, 1)^\mathrm{T}.
\end{equation*}
The normal cones are bounded by three hyperplanes whose normal vectors are
\begin{equation*}
n_{01} = (-1, 1, 0, 0, 1, -1)^\mathrm{T},
n_{12} = (0, -1, 1, 1, -1, 0)^\mathrm{T},
n_{20} = (1, 0, -1, -1, 0, 1)^\mathrm{T}.
\end{equation*}
\end{example}
}

\begin{example}
The normal vectors to the $22$ facets of $P_5$, and thus the rays of
the normal cones to $P_5$, form two classes (see Fig.
\ref{fig:p5rays}). The first class contains a total of $12$ vectors
(as there are $12$ assignments of nodes $0$ to $4$ to the labels $a$
to $e$ which yield nonisomorphic labelings), and every normal cone
contains six of them. The second class contains $10$ vectors, and
again every normal cone has six of them as rays.
\begin{figure}[htb]
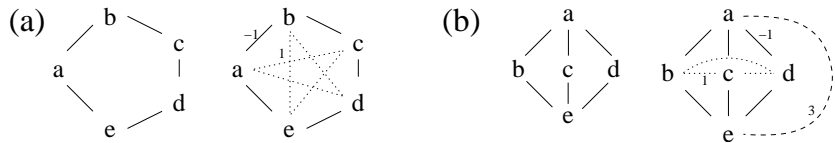

\vskip -0.1in
\begin{center}
\resizebox{0.9\textwidth}{!}{\input figs/p5rays.pstex_t}
\end{center}
\caption{Diagrams describing the facet-normals of $P_5$.}\label{fig:p5rays}
\vskip -0.2in
\end{figure}
For each class there are two diagrams in Fig. \ref{fig:p5rays}, and we
obtain a normal vector to one of the facets of $P_5$ by assigning
nodes from $\{0,\ldots,4\}$ to the labels $a,\ldots,e$. The left
diagram tells which vertices of $P_5$ belong to the facet defined by
that normal vector: Two nodes in the diagram are connected by an edge
iff the vertex belonging to that pair of nodes is in the facet. The
edges in the right diagram are labeled with the distance between the
corresponding pair of nodes in the normal vector.
This calls for an example: Setting $a = 0$, \dots, $e = 4$, Fig. \ref{fig:p5rays}(a) gives a distance vector
\[
(d_{01}, d_{02}, \ldots, d_{24}, d_{34})^\mathrm{T} = (-1, 1, 1, -1, -1, 1, 1, -1, 1, -1)^\mathrm{T},
\]
which is a common ray of the cones $cd_{01}$, $cd_{12}$, $cd_{23}$,
$cd_{34}$ and $cd_{04}$. Thus of the $22$ facets of $P_5$, $12$ have
five vertices and $10$ have six vertices.
\end{example}
\vskip -0.2in

\section{The NJ cones for five taxa}\label{C5}

\begin{figure}[ht]
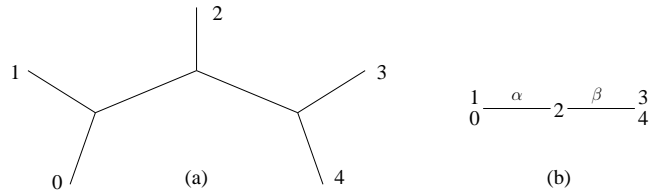

\vskip -0.2in
\begin{center}
\resizebox{0.7\textwidth}{!}{\input figs/fivetree.pstex_t}
\end{center}
\caption{(a) A tree with five taxa (b) The same tree with all edges
adjacent to leaves reduced to length zero. The remaining two edges
have lengths $\alpha$ and $\beta$.}
\label{fig:fivetree}
\vskip -0.2in
\end{figure}

In the case of five taxa there is just one unlabeled tree topology
(cf. Fig. \ref{fig:fivetree}) and there are 15 distinct labeled
trees: We have five choices for the leaf which is not part of a cherry
and then three choices how to group the remaining four leaves into two
pairs. For each of these labeled topologies, there are two ways in
which they might be reconstructed by the NJ algorithm: There are two
pairs, any one of which might be chosen in the first step of the NJ
algorithm.

For distinct leaf labels $a$, $b$ and $c \in \{0,1,2,3,4\}$ we define
$C_{ab,c}$ to be the set of all input vectors for which the cherry
$a$-$b$ is picked in the first step and $c$ remains as single node not
part of a cherry after the second step. For example, the tree in
Fig. \ref{fig:fivetree}(a) is the result for all vectors in
$C_{10,2} \cup C_{43,2}$. Since for each tree topology $((a, b), c, (d, e))$
(this tree topology is written in the Newick format) 
for distinct taxa $a, b, c, d, e \in \{0, 1, 2, 3, 4\}$, the NJ algorithm 
returns the same tree topology with any vector in the
union of two cones $C_{ab, c} \cup C_{de, c}$, there are 30 such cones in 
total, and we call the set of these cones $\mathcal{C}$.
\vskip -0.1in

\subsection{Permuting leaf labels}
\label{sec:leafperm}

Because there is only one unlabeled tree topology, we can map any
labeled topology to any other labeled topology by only changing the
labels of the leafs. Such a change of labels also permutes the entries
of the distance matrix. 
In this way, we get an action of the symmetric group $S_5$ on the input
space $R^{10}$, and the permutation $\sigma \in S_5$ maps the cone
$C_{ab,c}$ linearly to the cone
$C_{\sigma(a)\sigma(b),\sigma(c)}$. Therefore any property of the cone
$C_{ab,c}$ which is preserved by unitary linear transformations must be the
same for all cones in $\mathcal{C}$, and it suffices to determine it
for just one cone.

The action of $S_5$ on $R^{10}$ decomposes into irreducible representations by
\begin{equation*}
\underbrace{\R \oplus \R^4}_{= S} \oplus \underbrace{\R^5}_{=: W},
\end{equation*}
where the first summand is the subspace of all constant vectors and
the second one is the kernel of $A^{(5)}$. The sum of these two
subspaces is exactly the space $S$ generated by the shift vectors. The
third summand, which we call $W$, is the orthogonal complement of $S$
and it is spanned by vectors $w_{ab,cd}$ in $W$ with
\begin{equation*}
(w_{ab,cd})_{xy} := \begin{cases}
1 & \text{if }xy = ab\text{ or }xy = cd
\\
-1 & \text{if }xy = ac\text{ or }xy = bd
\\
0 & \text{else}
\end{cases}
\end{equation*}
where $a$, $b$, $c$ and $d$ are pairwise distinct taxa in $\{0, 1, 2, 3, 4\}$
and $(w_{ab,cd})_{xy}$ is the $x$--$y$th coordinate of the vector $w_{ab,cd}$. 
One linearly independent
subset of this is
\begin{equation*}
w_1 := w_{01,34},\quad
w_2 := w_{12,40},\quad
w_3 := w_{23,01},\quad
w_4 := w_{34,12},\quad
w_5 := w_{40,23}.
\end{equation*}
Note that the 5-cycle $(01234)$ of leaf labels cyclically permutes
these basis vectors, whereas the transposition $(01)$ acts via the
matrix
\begin{equation*}
T := \frac{1}{2}
\left(\begin{array}{rrrrr}
2 & 1 & 1 & 1 & 1 \\
0 & 1 &-1 &-1 &-1 \\
0 &-1 &-1 & 1 &-1 \\
0 &-1 & 1 &-1 &-1 \\
0 &-1 &-1 &-1 & 1
\end{array}\right).
\end{equation*}
Because a five-cycle and a transposition generate $S_5$, in principle
this gives us complete information about the operation.
\vskip -0.2in

\subsection{The cone $C_{43,2}$}

Since we can apply a permutation $\sigma \in S_5$, without loss of generality,
we suppose that the first cherry to be picked is cherry 9, which is the 
cherry with leaves 3 and 4. This is true for all input vectors 
$\mathbf{d}$ which satisfy
\begin{equation*}
(\mathbf{h}_{9,i}, \mathbf{d}) \geq 0
\text{ for }i=0,\ldots,8,
\end{equation*}
where the vector
\begin{equation*}
\mathbf{h}^{(n)}_{ij} := -A^{(n)}(\mathbf{e}_i - \mathbf{e}_j)
\end{equation*}
is perpendicular to the hyperplane of input vector for which cherries
$i$ and $j$ have the same Q-criterion, pointing into the direction of
vectors for which the Q-criterion of cherry $i$ is lower.

We let $\mathbf{r}_1$, $\mathbf{r}_2$ and $\mathbf{r}_3$ be the first
three rows of $-A^{(4)}R^{(5)}$. If $(\mathbf{r}_1, \mathbf{d})$ is
maximal then the second cherry to be picked is 0-1, leaving 2 as the
non-cherry node, and similarly $\mathbf{r}_2$ and $\mathbf{r}_3$ lead
to non-cherry nodes 1 and 0. This allows us to define the set of all
input vectors $\mathbf{d}$ for which the first picked cherry is 3-4
and the second one is 0-1:
\begin{equation}
\label{eqn:top34-2}
C_{34,2} := \{ \mathbf{d} \st
(\mathbf{h}_{9,i}, \mathbf{d}) \geq 0
\text{ for }i=0,\ldots,8,
\text{ and }
(\mathbf{r}_1-\mathbf{r}_2,\mathbf{d}) \geq 0,
(\mathbf{r}_1-\mathbf{r}_3,\mathbf{d}) \geq 0
\}.
\end{equation}

We have defined this set by 11 bounding hyperplanes. However, in fact, the
resulting cone has only nine facets. A computation using
\texttt{polymake} \cite{Gawrilow2000} reveals that the two hyperplanes 
$\mathbf{h}_{9,1}$
and $\mathbf{h}_{9,2}$ are no longer faces of the cone, while the
other nine hyperplanes in \eqref{eqn:top34-2} give exactly the facets
of the cone.  That means that, while we can find arbitrarily close
input vectors $\mathbf{d}$ and $\mathbf{d'}$ such that with an input
$\mathbf{d}$ the NJ algorithm will first pick cherry 3-4 and with an 
input $\mathbf{d'}$ the NJ algorithm will first pick cherry 1-2 
(or 0-2), we cannot do
this in such a way that $\mathbf{d}$ will result in the labeled tree
topology of Fig. \ref{fig:fivetree}, where 2 is the lonely leaf.

\comment{
\subsection{The Rays of $\mathcal{C}$}
\label{sec:raysofc}
Again using \texttt{polymake} \cite{Gawrilow2000}
we find that $C_{43,2}$ is the positive
cone spanned by fourteen rays. The union of the orbits of these rays
under the action of $S_5$ is a set of 82 vertices, which we call
$\mathcal{R}$. Each of the cones in $\mathcal{C}$ is spanned by a
certain subset with 14 vertices each of $\mathcal{R}$.
The set $\mathcal{R}$ has three orbits under the action of the
symmetric group $S_5$.  We characterize the rays in these orbits using the
graphs in Fig. \ref{fig:rays}. In these graphs, nodes which are
connected by an edge form a pair with the minimum Q-criterion in the first
step of the NJ algorithm, and the labels to each edge show which nodes 
are possible as the remaining unpaired node.
\begin{figure}[ht]
\vskip -0.2in
\begin{center}
\resizebox{0.9\textwidth}{!}{\input figs/rays.pstex_t}
\end{center}
\caption{Graphs characterizing the rays in $\mathcal{R}$. Here $a$,
$b$, $c$, $d$ and $e$ are pairwise distinct labels between 0 and 4.}
\label{fig:rays}
\vskip -0.2in
\end{figure}
For each graph in Fig. \ref{fig:rays} and each assignment of the
five leaves $0,\ldots,4$ to the variables $a,\ldots,e$ we get a set
$G$ of vectors which belong exactly to the cones indicated by the
graph. By what we have said in sections \ref{sec:shifting} and
\ref{sec:leafperm}, the set $G$ has the form
\begin{equation*}
  G = \{ \alpha\mathbf{g} + s \st \alpha \geq 0\text{ and }s\in S\},
\end{equation*}
where $S$ is the set of all shifting vectors in the sample space (which
is the subspace in $R^{10}$).
Thus $G$ is described by a vector $\mathbf{g}$ which is unique up to
normalization by a positive constant. We call this vector
$\mathbf{g}^{(\mathrm{a})}_{abcde}$,
$\mathbf{g}^{(\mathrm{b})}_{abcde}$ or
$\mathbf{g}^{(\mathrm{c})}_{abcde}$, depending on which graph it
corresponds to.
Some of these graphs coincide for different choices of
$a,\ldots,e$.
\begin{itemize}
\item The graphs in \ref{fig:rays}(a) are strings of four taxa. There
are five ways to choose the four taxa, and $\frac{4!}{2}=12$ ways in
which to arrange them (reversing the order of the string does not
matter). In total there are 60 rays of this type.
Thus for example the string 0-1-2-3 refers to an input vector
$\mathbf{g}^{(\mathrm{a})}_{01234}$ which
lies in the intersection of the five cones $C_{10,4}, C_{21,0}, C_{21,3}, 
C_{21,4}, C_{30, 4}$ but not in any other cone. 
It is worth noting that this vector is a
ray of the cones $C_{21,0}, C_{21,3}$ and $C_{21,4}$ but \emph{not} of 
the cones $C_{10,4}$
and $C_{32,4}$, which is indicated by the dotted edges in Fig.
\ref{fig:rays}. This implies that the cones in $\mathcal{C}$ do
\emph{not} form a fan, and in particular no polytope in $\R^{10}$ can have
$\mathcal{C}$ as its normal fan.
\item There are $\frac{4!}{2}$ possibilities for the cycle graph in
Fig. \ref{fig:rays}(b), each of which represents an input vector
which forms a ray of ten of the cones in $\mathcal{C}$.
\item Finally, there are $\binom{5}{2}=10$ graphs of type
\ref{fig:rays}(c), because these are determined by the two nodes we
choose for $a$ and $e$. Each of these corresponds to an input vector
that forms a ray of twelve cones.
\end{itemize}
If we collect all those rays in $\mathcal{R}$ which are rays of
$C_{43,2}$ we get
\begin{itemize}
\item 6 rays of type (a):
\begin{equation*}
\mathbf{g}^{(\mathrm{a})}_{03412},
\mathbf{g}^{(\mathrm{a})}_{03421},
\mathbf{g}^{(\mathrm{a})}_{13402},
\mathbf{g}^{(\mathrm{a})}_{13420},
\mathbf{g}^{(\mathrm{a})}_{23401},
\mathbf{g}^{(\mathrm{a})}_{23410}
\end{equation*}
\item 4 of type (b):
\begin{equation*}
\mathbf{g}^{(\mathrm{b})}_{34201},
\mathbf{g}^{(\mathrm{b})}_{34210},
\mathbf{g}^{(\mathrm{b})}_{34012},
\mathbf{g}^{(\mathrm{b})}_{34102}
\end{equation*}
\item 4 of type (c):
\begin{equation*}
\mathbf{g}^{(\mathrm{c})}_{34201},
\mathbf{g}^{(\mathrm{c})}_{34210},
\mathbf{g}^{(\mathrm{c})}_{34021},
\mathbf{g}^{(\mathrm{c})}_{34120}.
\end{equation*}
\end{itemize}
We call these vectors $\mathbf{g}_1,\ldots,\mathbf{g}_{14}$.
Any vector $\mathbf{v} \in C_{43,2}$ is a linear combination of these
fourteen vectors with nonnegative linear coefficients, say
\begin{equation*}
\mathbf{v} = \sum_i \alpha_i \mathbf{g}_i,
\end{equation*}
where we number the fourteen rays arbitrarily. By linearity, in both
steps of the NJ algorithm, the Q-criteria are linear combinations of 
the Q-criteria
for the $\mathbf{g}_i$, using the same linear coefficients
$\alpha_i$. Furthermore, because all of the vectors $\mathbf{g}_i$ lie
in $C_{43,2}$, in both steps of the NJ algorithm there is a pair of leaves 
such that
the minimum of the Q-criteria is attained at this pair for all
$\mathbf{g}_i$, and therefore the minimum Q-criterion of the linear
combination $\sum \alpha_i \mathbf{g}_i$ is the linear combination of
the minimum Q-criteria of the individual $\mathbf{g}_i$. This minimum
is attained exactly by those pairs which have minimum Q-criterion in all
$\mathbf{g}_i$ with strictly positive linear coefficient $\alpha_i$,
we get:
\begin{equation*}
\mathbf{v} \in C_{ab,c}
\quad
\text{iff}
\quad
\mathbf{g}_i \in C_{ab,c}\text{ for all }i\text{ with }\alpha_i > 0.
\end{equation*}
Therefore the rays in $\mathcal{R}$ are the input vectors for which the NJ
algorithm is least stable. We can give explicit descriptions of these vectors
using the graphs in Fig. \ref{fig:rayvectors}. The graphs give the
distance value assigned to each pair of leaves. For example, we get
\begin{equation*}
\mathbf{g}^{(\mathrm{a})}_{01234} = (-3, 5, -3, -1, 5, -3, -1, 1, 1, -1 )^\mathrm{T}.
\end{equation*}
\begin{figure}[ht]
\begin{center}
\resizebox{0.9\textwidth}{!}{\input figs/rayvectors.pstex_t}
\end{center}
\caption{Graphs describing the vectors $\mathbf{g}_{abcde}$ for each
of the three types of graphs. Labels assigned to edges denote the
distance assigned to the corresponding pair of leaves in the vector,
and edges of the same type receive the same label.}
\label{fig:rayvectors}
\vskip -0.2in
\end{figure}
\subsection{The Combinatorics of $C_{ab,c}$}
Each of the cones $C_{ab,c}$ is a five-dimensional cone, and
intersecting it with a suitable hyperplane leaves us with a
four-dimensional polytope $P$. This polytope has 14 vertices, one for
each ray of $C_{ab,c}$. The $F$-vector of the polytope is
$(14,32,27,9)$, so there are 32 edges, 27 ridges and 9 facets. For
details about polyhedral geometry see \cite{Ziegler95}.
\begin{figure}
\begin{center}
\resizebox{0.9\textwidth}{!}{\input figs/facetdiag1.pstex_t}
\end{center}
\caption{(a) and (b): The common facet $F_1$ of $C_{ab,c}$ and
$C_{de,c}$ and its graph. (c) and (d): The facets $F_2$ and $F_3$ and
their graph.}
\label{fig:facetdiag1}
\vskip -0.2in
\end{figure}
In section \ref{sec:raysofc} we gave a detailed description of the
rays of $C_{ab,c}$, which we have divided into three types. We shall
now do the same for the nine facets of $P$.
\begin{itemize}
\item The eight vertices of types (b) and (c) form a facet $F_1$,
which is the common facet of $C_{ab,c}$ and $C_{de,c}$, joining the
two cones resulting in the same topology. The graph of this facet and
a combinatorially equivalent realization of it are shown in Fig.
\ref{fig:facetdiag1}. It looks like a slightly distorted
three-dimensional cube with eight facets. In fact, each of these eight
facets is shared with exactly one other facet of $P$, making $F_1$ a
good starting-point for a Schlegel-diagram of $P$.
\item Two of the square facets of $F_1$ belong to the common facets of
$C_{ab,c}$ with $C_{ab,d}$ and $C_{ab,e}$ respectively. We call these
$F_2$ and $F_3$. In turn, these two share a common facet, which
consists of the six type-(a)-rays of $C_{ab,c}$ joined in a hexagon.
\item The other two square facets of $F_1$ form triangular prisms $F_4$ and $F_5$
with two of the type-(a)-rays each.
\item The remaining four triangular facets of $F_1$ form egyptian
pyramids $F_6$ to $F_9$ with two of the type-(a)-rays each.
\end{itemize}
Using $F_1$ as a base, we get the Schlegel-diagram shown in Fig. \ref{fig:schlegel}.
\begin{figure}
\begin{center}
\resizebox{0.5\textwidth}{!}{\input figs/schlegel.pstex_t}
\end{center}
\caption{The Schlegel-diagram of $P$.}
\label{fig:schlegel}
\vskip -0.2in
\end{figure}
}

\comment{
\subsection{The Cone of Tree Metrics}
Suppose our input vector is a tree metric. By relabeling the tree and
using the shifting lemma, we can reduce the tree to Fig.
\ref{fig:fivetree}(b), where all edges connecting to leaves have
length zero. The edge-lengths $\alpha$ and $\beta$ determine the
distance vector $\mathbf{v}$ as follows:
\begin{equation*}
\mathbf{v} = \alpha(0, 1, 1, 1, 1, 0, 1, 1, 0, 0)^\mathrm{T}
+
\beta(0, 0, 0, 1, 1, 1, 1, 1, 1, 0)^\mathrm{T}
=: \alpha\mathbf{a} + \beta\mathbf{b}.
\end{equation*}
If we assume these to be non-negative, we get a cone spanned by the
two vectors $\mathbf{a}$ and $\mathbf{b}$. These satisfy
\begin{equation*}
\mathbf{a} \in C_{10,2} \cap C_{10,3} \cap C_{10,4},
\quad
\mathbf{b} \in C_{43,0} \cap C_{43,1} \cap C_{43,2}
\quad\text{and}\quad
\underbrace{\frac{1}{2}(\mathbf{a}+\mathbf{b})}_{=: \mathbf{b}'} \in C_{10,2} \cap C_{43,2},
\end{equation*}
and none of these vectors lies in any other cone.
We normalize to $\alpha+\beta = 1$, and because of symmetry we may
assume that $\alpha \geq \beta$, or equivalently $\alpha \geq
\frac{1}{2}$. Then we can write $\mathbf{v}$ as
\begin{align*}
\mathbf{v} 
&= \alpha \mathbf{a} + (1-\alpha)\mathbf{b}
\\
&= \left(\alpha-\frac{1}{2}\right) \mathbf{a} + (1-\alpha)\mathbf{b}'
\end{align*}
Because both $\mathbf{a}$ and $\mathbf{b}'$ are in the convex set
$C_{10,2}$, the convex linear combination $\mathbf{v}$ is also in
$C_{10,2}$, proving correctness of the NJ algorithm for five taxa.
We can also compute the $\ell_2$ distance of $\mathbf{v}$ from the
faces of the cone $C_{10,2}$. There are 9 hyperplanes defining $C_{10,2}$, 
but we may ignore one of them, which defines the common facet with
$C_{43,2}$. For any hyperplane $H$ such that $\mathbf{a}$ and
$\mathbf{b}'$ lie on the same side of $H$, the distance
$d(\mathbf{v},H)$ between $H$ and $\mathbf{v}$ is given by
\begin{equation*}
d(\mathbf{v},H) = \left(\alpha-\frac{1}{2}\right) d(\mathbf{a},H) + (1-\alpha)d(\mathbf{b}',H),
\end{equation*}
and taking the minimum of this over the eight remaining faces of
$C_{10,2}$ we obtain
\begin{equation*}
d(\mathbf{v},(C_{10,2}\cap C_{43,2})^\mathrm{c}) = \frac{1-\alpha}{\sqrt{3}},
\end{equation*}
where $(C_{10,2}\cap C_{43,2})^\mathrm{c}$ is the complement of 
$C_{10,2}\cap C_{43,2}$.
If we divide this by the length $\beta = 1-\alpha$ of the smaller of
the two interior edges, we get an $\ell_2$-radius of $1/\sqrt{3}
\approx 0.577$. Note that because our method relies on orthogonal
projections, we get $\ell_2$ bounds instead of $\ell_\infty$ bounds.
}

Also note that $\mathcal{C}$, the set of NJ cones which the NJ algorithm 
returns the same tree topology with any vector in the
union of two cones $C_{ab, c} \cup C_{de, c}$, is not convex which is shown
in \cite{NJME}.  For details of geometry and combinatorics of the NJ 
cones for $n = 5$, see \cite{NJME}.

\vskip -0.2in

\section{The six taxa case}\label{6taxa}

\begin{figure}
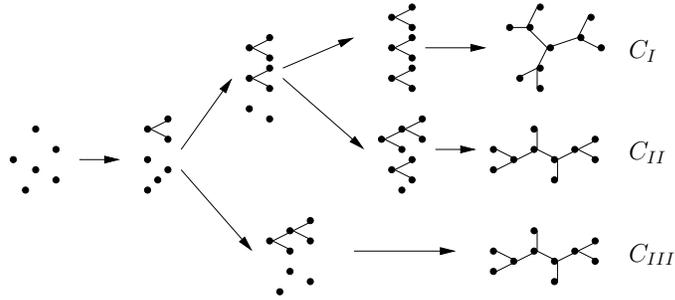

\begin{center}
\resizebox{0.8\textwidth}{!}{\input figs/sixcones.pstex_t}
\end{center}
\caption{The three ways of picking cherries in the six taxa case.}
\label{fig:sixcones}
\vskip -0.1in
\end{figure}

Note that since each of the NJ cones includes constraints for five taxa,
the union of the NJ cones which gives the same tree topology is not
convex.  
To analyze the behavior of NJ on distance matrices for six taxa, we
use the action of the symmetric group as much as
possible. However, in this case we get three different classes of
cones which cannot be mapped onto each other by this action. We assume
the cherry which is picked in the first step to consist of the nodes
$4$ and $5$. Picking this cherry replaces these two nodes by a newly
created node $45$, and we have to distinguish two different cases in
the second step (see Fig. \ref{fig:sixcones}):
\begin{itemize}
\item If the cherry in the second step does not contain the new node
$45$, we may assume the cherry to be $01$. For the third step, we
again get two possibilities:
\begin{itemize}
\item The two nodes $45$ and $01$ get joined in the third step. We
call the cone of input vectors for which this happens $C_\mathrm{I}$.
\item The node $45$ is joined to one of the nodes $2$ and $3$, without loss of
generality, to $3$. We call the resulting cone $C_\mathrm{II}$.
\end{itemize}
\item If the cherry in the second step contains the new node $45$, we
may assume the other node of this cherry to be $3$, creating a new
node $45-3$. In the third step, all that matters is which of the three
nodes $0$, $1$ and $2$ is joined to the node $45-3$, and we may, without loss of
generality,
assume this to be node $2$. This gives the third type of cone,
$C_\mathrm{III}$.
\end{itemize}

\begin{figure}
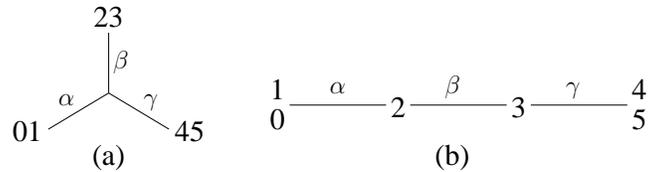

\begin{center}
\resizebox{0.7\textwidth}{!}{\input figs/sixtree.pstex_t}
\end{center}
\caption{The two possible topologies for trees with six leaves, with
edges connecting to leaves shrunk to zero.}
\label{fig:sixtree}
\vskip -0.2in
\end{figure}

The resulting tree topology for the cone $C_\mathrm{I}$ is shown in
Fig. \ref{fig:sixtree}(a), while both $C_\mathrm{II}$ and
$C_\mathrm{III}$ give the topology shown in \ref{fig:sixtree}(b). We
now determine which elements of $S_6$ leave these cones fixed
(stabilizer) and how many copies of each cone give the same labeled
tree topology:

\begin{center}
\begin{tabular}{l|c|c|c}
&$C_\mathrm{I}$&$C_\mathrm{II}$&$C_\mathrm{III}$
\\
\hline
stabilizer & $\langle (01), (23), (45) \rangle$ & $\langle (01),(45) \rangle$ & $\langle (01),(45) \rangle$
\\
size of stabilizer & 8 & 4 & 4
\\
number of cones & 90 & 180 & 180
\\
cones giving same labeled topology
& 6 & 2 & 2
\\
solid angle (approx.)
& $2.888\cdot 10^{-3}$
& $1.848\cdot 10^{-3}$
& $2.266\cdot 10^{-3}$
\end{tabular}
\end{center}

Thus, the input space $\R^{15}$ is divided into 450 cones, 90 of type
I and 180 each of types II and III. There are $15$ different ways of
assigning labels to the tree topology in Fig. \ref{fig:sixtree}(a),
and for each of these there are six copies of $C_\mathrm{I}$ whose
union describes the set of input vectors resulting in that
topology. For the topology in Fig. \ref{fig:sixtree}(b) we get 90
ways of assigning labels to the leaves, each corresponding to a union
of two copies of $C_\mathrm{II}$ and two copies of $C_\mathrm{III}$.

The above table also gives the solid angles of the three cones. In the
five taxa case, any two cones can be mapped onto one another by the
action of the symmetric group, which is unitary. Therefore all thirty
cones have the same solid angle, which must be $1/30$. However, in the
six taxa case, we get different solid angles, and we see that about
three thirds of the solid angle at the origin are taken by the cones
of types $II$ and $III$. Thus, on a random vector chosen according to
any probability law which is symmetric around the origin (e.g.,
standard normal distribution), NJ will output the tree topology of
Fig. \ref{fig:sixtree}(b) with probability about $3/4$.

On the other hand, any labeled topology of the type in Fig.
\ref{fig:sixtree}(a) belongs to six cones of type $I$, giving a total
solid angle of $\approx 1.73\cdot 10^{-2}$, whereas any labeled
topology of the type in Fig. \ref{fig:sixtree}(b) belongs to two
cones each of type $II$ and $III$, giving a total solid angle of only
$\approx 0.82 \cdot 10^{-2}$, which is half as much. This suggests
that reconstructing trees of the latter topology is less robust against
noisy input data.
\vskip -0.2in

\section{Simulation results}\label{simulations}
\begin{figure}[ht]
\vskip -0.3in
\begin{center}
\includegraphics[scale=0.4]{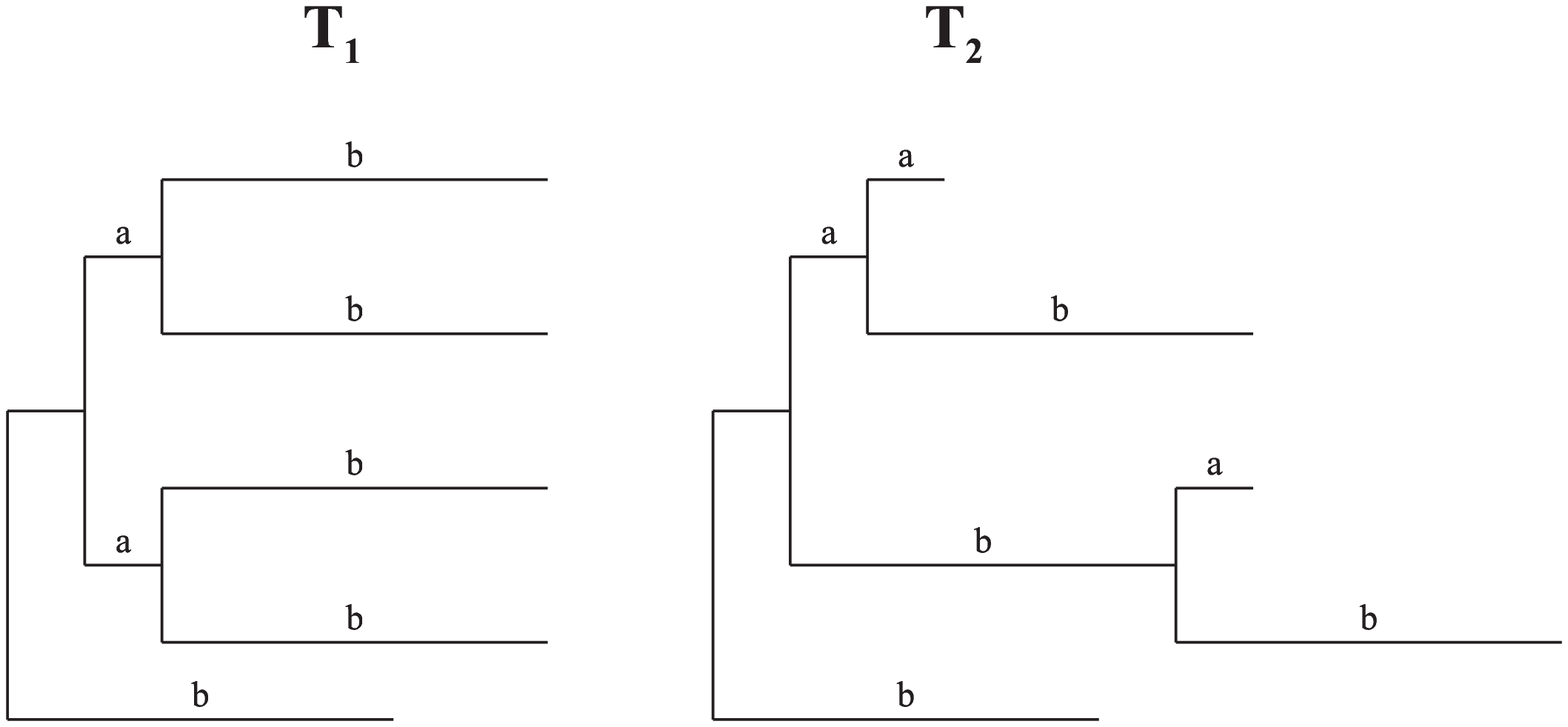}
\end{center}
\caption{$T_1$ and $T_2$ tree models which are subtrees of the tree
models in \cite{Ota2000}.}
\label{fig:T1T2}
\vskip -0.2in
\end{figure}

In this section we will analyze how the tree metric for a tree and
pairwise distances estimated via the maximum likelihood estimation
lie in the polyhedral subdivision of the sample space.  
In particular, we analyze
subtrees of the two parameter family of trees described by
\cite{Ota2000}.  These are trees for which the NJ algorithm has
difficulty in resolving the correct topology.  In order to understand
which cones the data lies in, we simulated 10,000 data sets on each
of the two tree shapes, $T_1$ and $T_2$ (Fig. \ref{fig:T1T2}) at the
edge length ratio, a/b = 0.03/0.42 for sequences of length 500BP under
the Jukes-Cantor model \cite{Jukes1969}. We also repeated the runs
with the Kimura 2-parameter model \cite{Kimura1980}. They are the
cases (on eight taxa) in \cite{Ota2000} that the NJ algorithm had
most difficulties in their simulation study (also the same as in
\cite{Levy2005}).  Each set of 5 sequences are generated via {\tt
evolver} from {\tt PAML} package \cite{Yang1997} under the given
model.  {\tt evolver} generates a set of sequences given the model and
tree topology using the birth-and-death process.  For each set of 5
sequences, we compute first pairwise distances via the heuristic MLE
method using a software {\tt fastDNAml} \cite{Olsen1994}.
To compute cones, we used {\tt MAPLE} and \texttt{polymake}.

\begin{figure}[!ht]
\resizebox{0.9\textwidth}{!}{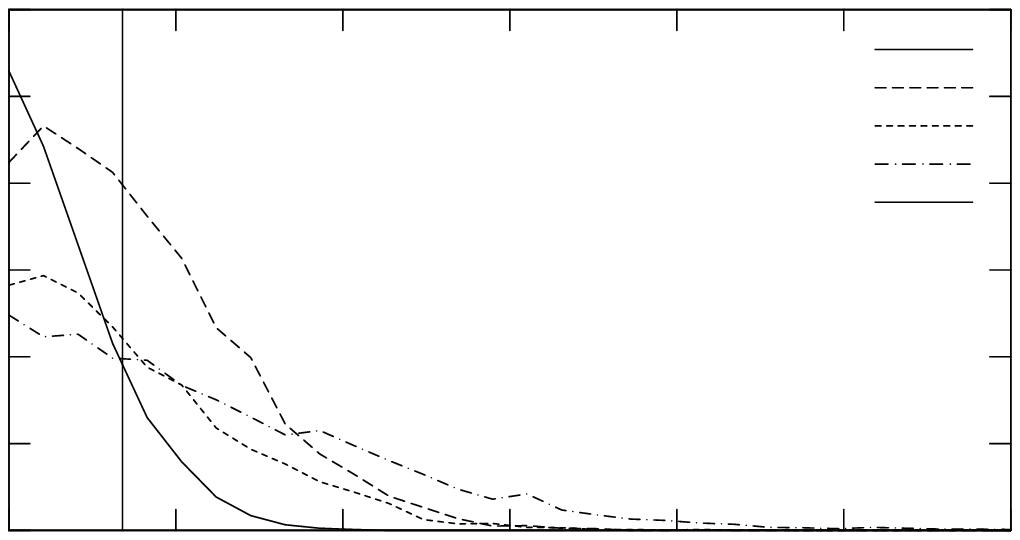}\\
\resizebox{0.9\textwidth}{!}{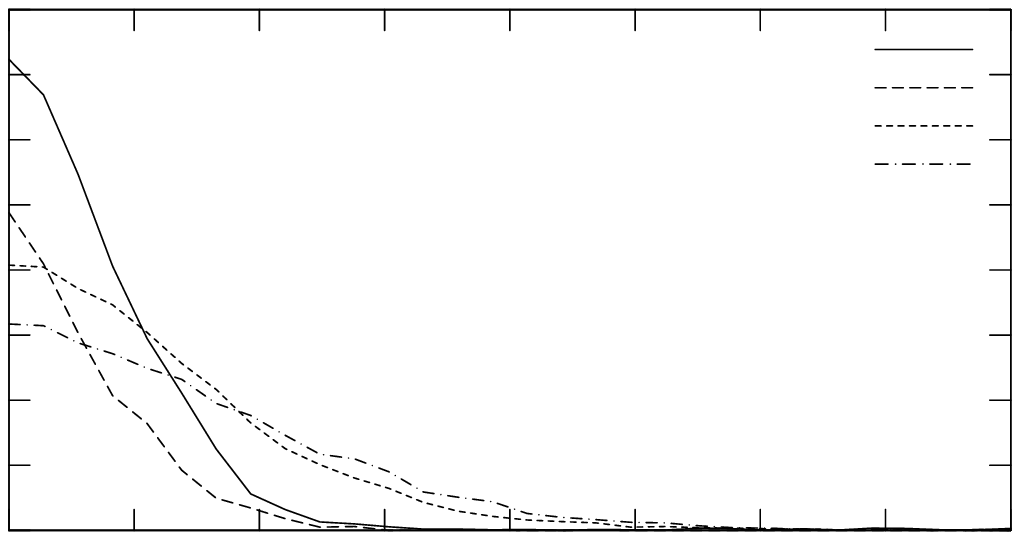}
\caption{Distances of correctly (top) and incorrectly (bottom)
classified input vectors to the closest incorrectly/correctly
classified vector.}
\label{fig:simplots}
\vskip -0.2in
\end{figure}

To study how far each set of pairwise distances estimated via the
maximum likelihood estimation (which is a vector $\mathbf{y}$ in
$\R^5$) lies from the cone, where the additive tree metric lies, in
the sample space, we calculated the $\ell_2$-distance between the cone
and a vector $\mathbf{y}$.
\begin{figure}[ht]
\vskip -0.2in
\begin{center}
\begin{tabular}{l|cc|cc}
&\multicolumn{2}{c}{\textbf{JC}}&\multicolumn{2}{c}{\textbf{Kimura2}}
\\
&T1&T2&T1&T2
\\
\hline
\textbf{\# of cases}& 3,581 & 6,441 & 3,795 & 4,467
\\
\textbf{Mean}&$0.0221$&$0.0421$ & $0.0415$ & $0.0629$
\\
\textbf{Variance}& $2.996\cdot 10^{-4}$ & $9.032\cdot 10^{-4} $
&
$1.034 \cdot 10^{-3}$ & $2.471 \cdot 10^{-3}$
\end{tabular}
\end{center}
\caption{Mean and variance of the
distances of correctly classified
vectors from the nearest
misclassified vector.}
\label{fig:goodstat}
\end{figure}

\begin{figure}[!ht]
\begin{center}
\begin{tabular}{l|cc|cc}
&\multicolumn{2}{c}{\textbf{JC}}&\multicolumn{2}{c}{\textbf{Kimura2}}
\\
&T1&T2&T1&T2
\\
\hline
\textbf{\# of cases}& 6,419 & 3,559 & 6,205 & 5,533
\\
\textbf{Mean}&$0.0594$&$0.0331$ & $0.0951$ & $0.0761$
\\
\textbf{Variance}& $0.0203$ & $7.39\cdot 10^{-4} $
&
$0.0411$ & $3.481 \cdot 10^{-3}$
\end{tabular}
\end{center}
\caption{Mean and variance of the
distances of misclassified vectors to the nearest correctly classified vector.}
\label{fig:badstat}
\vskip -0.2in
\end{figure}

Suppose we have a cone $C$ defined by hyperplanes
$\mathbf{n}_1,\ldots,\mathbf{n}_r$, i.e.,
\begin{equation*}
C = \{ \mathbf{x} \st (\mathbf{n}_i,\mathbf{x}) \geq 0\text{ for
}i=1,\ldots,r\},
\end{equation*}
and we want to find the closest point in $C$ from some given point
$\mathbf{v}$. Because $C$ is convex, for $\ell_2$-norm there is only
one such point, which we call $\mathbf{u}$. If $\mathbf{v} \in C$ then
$\mathbf{u} = \mathbf{v}$ and we are done. If not, there is at least
one $\mathbf{n}_i$ with $(\mathbf{n}_i, \mathbf{v}) < 0$, and
$\mathbf{u}$ must satisfy $(\mathbf{n}_i,\mathbf{u}) = 0$.

Now the problem reduces to a lower dimensional problem of the same
kind: We project $\mathbf{v}$ orthogonally into the hyperplane $H$
defined by $(\mathbf{n}_i,\mathbf{x}) = 0$ and call the new vector
$\tilde{\mathbf{v}}$. Also, $C \cap H$ is a facet of $C$, and in
particular a cone, so we proceed by finding the closest point in this
cone from $\tilde{\mathbf{v}}$.

We say an input vector (distance matrix) is {\em correctly classified}
if the vector is in one of the cones where the vector representation 
of the tree metric (noiseless input) is.
We say an input vector is {\em incorrectly classified}
if the vector is in the complement of the cones where 
the vector representation of the tree metric is.
For input vectors (distance matrices) which are correctly classified by 
the NJ algorithm, we compute the
minimum distance to any cone giving a different tree topology. This
distance gives a measure of robustness or confidence in the result,
with bigger distances meaning greater reliability. The results are
plotted in the left half of Fig. \ref{fig:simplots} and in Fig.
\ref{fig:goodstat}. Note that the distance of the noiseless input,
i.e., the tree metric from the tree we used for
generating the data samples, gives an indication of what order of
magnitude to expect with these values.

For input vectors to which the NJ algorithm returns with a tree topology 
different
from the correct tree topology, we compute the
distances to the two cones for which the correct answer is given and
take the minimum of the two. The bigger this distance is, the further
we are off. The results are shown in the right half of Fig.
\ref{fig:simplots} and in Fig. \ref{fig:badstat}.

From our results in Fig. \ref{fig:goodstat} and Fig. \ref{fig:badstat},
one notices that the NJ algorithm returns the correct tree more often with 
$T_2$ than with $T_1$.  These results are consistent with the results in 
\cite{Steel2006,Mihaescu2006}.  Note that any possible quartet in $T_1$
has a smaller (or equal) length of its internal edge than in $T_2$ 
(see Fig. \ref{fig:T1T2}).  Gascuel and Steel defined this measure 
as {\em neighborliness} \cite{Steel2006}.  Mihaescu et al. showed that the NJ algorithm
returns the correct tree if it works correctly locally for the quartets in
the tree \cite{Mihaescu2006}.  The neighborliness of a quartet is one of the most important
factors to reconstruct the quartet correctly, i.e., the shorter it is the more 
difficult the NJ algorithm returns the correct quartet. 
Also Fig. \ref{fig:simplots} shows that most of the input vectors lie around
the  boundary of cones, including the noiseless input vector (the tree 
metric).  This shows that the tree models $T_1$ and $T_2$ are difficult
for the NJ algorithms to reconstruct the correct trees.  
All source code for these simulations described in this paper 
will be available at authors' websites.
\vskip -0.2in

\section{Open problems}\label{problem}
\comment{
There has been much research on statistical
consistency and the robustness of the NJ algorithm to small perturbations
of tree metrics, such as \cite{Studier1988,Atteson99,Mihaescu2006}. 
This geometric view of the NJ algorithm provides new perspectives
of the questions on statistical
consistency and the robustness of the NJ algorithm to small perturbations
of tree metrics.
We conclude with an open problem
related to the topics in this paper.}
\vskip -0.1in

\begin{question}
Can we use the NJ cones for analyzing how the NJ algorithm works if each pairwise distance is assumed
to be of the form $D_0 + \epsilon$ where $D_0$ is the unknown true tree 
metric, and $\epsilon$ is a collection of independent normally distributed
random variables? 
We think this would be very interesting and relevant.
\end{question}
\begin{question}
With any $n$, is there an efficient method for computing (or approximating) 
the distance between a given pairwise distance vector and the boundary of the NJ 
optimality region which contains it?  This problem is equivalent to projecting
a point inside a {\em polytopal complex} $P$ onto the boundary of $P$.
Note that the size of the complex grows very fast with $n$.  How fast
does the number of the complex grow?
This would allow assigning a confidence score to the tree topology computed
by the NJ algorithm.
\end{question}

%
%

%
\end{document}

%% file: figs/sampletree.pstex_t
\begin{picture}(0,0)%
\includegraphics{sampletree.pstex}%
\end{picture}%
\setlength{\unitlength}{4144sp}%
\begingroup\makeatletter\ifx\SetFigFont\undefined%
\gdef\SetFigFont#1#2#3#4#5{%
  \reset@font\fontsize{#1}{#2pt}%
  \fontfamily{#3}\fontseries{#4}\fontshape{#5}%
  \selectfont}%
\fi\endgroup%
\begin{picture}(1745,781)(6103,-3939)
\put(6211,-3256){\makebox(0,0)[rb]{\smash{{\SetFigFont{14}{16.8}{\familydefault}{\mddefault}{\updefault}{\color[rgb]{0,0,0}a}%
}}}}
\put(6211,-3931){\makebox(0,0)[rb]{\smash{{\SetFigFont{14}{16.8}{\familydefault}{\mddefault}{\updefault}{\color[rgb]{0,0,0}b}%
}}}}
\put(7741,-3256){\makebox(0,0)[lb]{\smash{{\SetFigFont{14}{16.8}{\familydefault}{\mddefault}{\updefault}{\color[rgb]{0,0,0}c}%
}}}}
\put(7741,-3931){\makebox(0,0)[lb]{\smash{{\SetFigFont{14}{16.8}{\familydefault}{\mddefault}{\updefault}{\color[rgb]{0,0,0}d}%
}}}}
\put(6526,-3346){\makebox(0,0)[lb]{\smash{{\SetFigFont{11}{13.2}{\familydefault}{\mddefault}{\updefault}{\color[rgb]{0,0,0}1}%
}}}}
\put(6526,-3886){\makebox(0,0)[lb]{\smash{{\SetFigFont{11}{13.2}{\familydefault}{\mddefault}{\updefault}{\color[rgb]{0,0,0}2}%
}}}}
\put(7471,-3886){\makebox(0,0)[rb]{\smash{{\SetFigFont{11}{13.2}{\familydefault}{\mddefault}{\updefault}{\color[rgb]{0,0,0}1}%
}}}}
\put(7426,-3346){\makebox(0,0)[rb]{\smash{{\SetFigFont{11}{13.2}{\familydefault}{\mddefault}{\updefault}{\color[rgb]{0,0,0}0.3}%
}}}}
\put(6976,-3481){\makebox(0,0)[b]{\smash{{\SetFigFont{11}{13.2}{\familydefault}{\mddefault}{\updefault}{\color[rgb]{0,0,0}0.5}%
}}}}
\end{picture}%

%% file: figs/p5rays.pstex_t
\begin{picture}(0,0)%
\includegraphics{p5rays.pstex}%
\end{picture}%
\setlength{\unitlength}{4144sp}%
\begingroup\makeatletter\ifx\SetFigFont\undefined%
\gdef\SetFigFont#1#2#3#4#5{%
  \reset@font\fontsize{#1}{#2pt}%
  \fontfamily{#3}\fontseries{#4}\fontshape{#5}%
  \selectfont}%
\fi\endgroup%
\begin{picture}(9753,1581)(940,-3976)
\end{picture}%

%% file: figs/fivetree.pstex_t
\begin{picture}(0,0)%
\includegraphics{fivetree.pstex}%
\end{picture}%
\setlength{\unitlength}{4144sp}%
\begingroup\makeatletter\ifx\SetFigFont\undefined%
\gdef\SetFigFont#1#2#3#4#5{%
  \reset@font\fontsize{#1}{#2pt}%
  \fontfamily{#3}\fontseries{#4}\fontshape{#5}%
  \selectfont}%
\fi\endgroup%
\begin{picture}(6858,2013)(2939,-3199)
\put(9226,-2221){\makebox(0,0)[b]{\smash{{\SetFigFont{12}{14.4}{\familydefault}{\mddefault}{\updefault}{\color[rgb]{0,0,0}$\beta$}%
}}}}
\put(8371,-2221){\makebox(0,0)[b]{\smash{{\SetFigFont{12}{14.4}{\familydefault}{\mddefault}{\updefault}{\color[rgb]{0,0,0}$\alpha$}%
}}}}
\end{picture}%

%% file: figs/rays.pstex_t
\begin{picture}(0,0)%
\includegraphics{rays.pstex}%
\end{picture}%
\setlength{\unitlength}{4144sp}%
\begingroup\makeatletter\ifx\SetFigFont\undefined%
\gdef\SetFigFont#1#2#3#4#5{%
  \reset@font\fontsize{#1}{#2pt}%
  \fontfamily{#3}\fontseries{#4}\fontshape{#5}%
  \selectfont}%
\fi\endgroup%
\begin{picture}(10605,1624)(-1777,-3946)
\end{picture}%

%% file: figs/rayvectors.pstex_t
\begin{picture}(0,0)%
\includegraphics{rayvectors.pstex}%
\end{picture}%
\setlength{\unitlength}{4144sp}%
\begingroup\makeatletter\ifx\SetFigFont\undefined%
\gdef\SetFigFont#1#2#3#4#5{%
  \reset@font\fontsize{#1}{#2pt}%
  \fontfamily{#3}\fontseries{#4}\fontshape{#5}%
  \selectfont}%
\fi\endgroup%
\begin{picture}(10990,1611)(-1737,-3946)
\end{picture}%

%% file: figs/facetdiag1.pstex_t
\begin{picture}(0,0)%
\includegraphics{facetdiag1.pstex}%
\end{picture}%
\setlength{\unitlength}{4144sp}%
\begingroup\makeatletter\ifx\SetFigFont\undefined%
\gdef\SetFigFont#1#2#3#4#5{%
  \reset@font\fontsize{#1}{#2pt}%
  \fontfamily{#3}\fontseries{#4}\fontshape{#5}%
  \selectfont}%
\fi\endgroup%
\begin{picture}(6404,1470)(3599,-3841)
\end{picture}%

%% file: figs/sixcones.pstex_t
\begin{picture}(0,0)%
\includegraphics{sixcones.pstex}%
\end{picture}%
\setlength{\unitlength}{4144sp}%
\begingroup\makeatletter\ifx\SetFigFont\undefined%
\gdef\SetFigFont#1#2#3#4#5{%
  \reset@font\fontsize{#1}{#2pt}%
  \fontfamily{#3}\fontseries{#4}\fontshape{#5}%
  \selectfont}%
\fi\endgroup%
\begin{picture}(6495,2590)(4691,-4866)
\put(10171,-2761){\makebox(0,0)[lb]{\smash{{\SetFigFont{14}{16.8}{\familydefault}{\mddefault}{\updefault}{\color[rgb]{0,0,0}$C_{I}$}%
}}}}
\put(10171,-3661){\makebox(0,0)[lb]{\smash{{\SetFigFont{14}{16.8}{\familydefault}{\mddefault}{\updefault}{\color[rgb]{0,0,0}$C_{II}$}%
}}}}
\put(10171,-4561){\makebox(0,0)[lb]{\smash{{\SetFigFont{14}{16.8}{\familydefault}{\mddefault}{\updefault}{\color[rgb]{0,0,0}$C_{III}$}%
}}}}
\end{picture}%

%% file: figs/sixtree.pstex_t
\begin{picture}(0,0)%
\includegraphics{sixtree.pstex}%
\end{picture}%
\setlength{\unitlength}{4144sp}%
\begingroup\makeatletter\ifx\SetFigFont\undefined%
\gdef\SetFigFont#1#2#3#4#5{%
  \reset@font\fontsize{#1}{#2pt}%
  \fontfamily{#3}\fontseries{#4}\fontshape{#5}%
  \selectfont}%
\fi\endgroup%
\begin{picture}(4793,1293)(4734,-2839)
\put(8056,-2221){\makebox(0,0)[b]{\smash{{\SetFigFont{12}{14.4}{\familydefault}{\mddefault}{\updefault}{\color[rgb]{0,0,0}$\beta$}%
}}}}
\put(7201,-2221){\makebox(0,0)[b]{\smash{{\SetFigFont{12}{14.4}{\familydefault}{\mddefault}{\updefault}{\color[rgb]{0,0,0}$\alpha$}%
}}}}
\put(8956,-2221){\makebox(0,0)[b]{\smash{{\SetFigFont{12}{14.4}{\familydefault}{\mddefault}{\updefault}{\color[rgb]{0,0,0}$\gamma$}%
}}}}
\put(5176,-2311){\makebox(0,0)[b]{\smash{{\SetFigFont{12}{14.4}{\familydefault}{\mddefault}{\updefault}{\color[rgb]{0,0,0}$\alpha$}%
}}}}
\put(5581,-2041){\makebox(0,0)[b]{\smash{{\SetFigFont{12}{14.4}{\familydefault}{\mddefault}{\updefault}{\color[rgb]{0,0,0}$\beta$}%
}}}}
\put(5806,-2311){\makebox(0,0)[b]{\smash{{\SetFigFont{12}{14.4}{\familydefault}{\mddefault}{\updefault}{\color[rgb]{0,0,0}$\gamma$}%
}}}}
\end{picture}%

%% file: AB2008_v5.bbl
\begin{thebibliography}{5}
%
\bibitem{Atteson99}
Atteson, K.:
The performance of neighbor-joining methods of phylogenetic reconstruction.
{Algorithmica}, {\bf 25} (1999) 251--278.
\bibitem{Bryant2005}
Bryant, D.:
On the uniqueness of the selection criterion in neighbor-joining.
J. Classif. {\bf 22} (2005) 3--15.

\bibitem{NJME}
Eickmeyer, K., Huggins, P., Pachter, L. and Yoshida, R.:
On the optimality of the neighbor-joining algorithm.
To appear in Algorithms in Molecular Biology. 

\bibitem {Felsenstein1981}
Felsenstein, J.:
Evolutionary trees from {DNA} sequences: a maximum likelihood approach.
Journal of Molecular Evolution {\bf 17} (1981) 368--376.
%
\bibitem {Galtier2005}
Galtier, N., Gascuel, O., and Jean-Marie, A.:
Markov models in molecular evolution.
In {\em Statistical Methods in Molecular Evolution} edited by Nielsen, R.,
(2005) 3--24.
\bibitem{Gawrilow2000}
Gawrilow, E. and Joswig, M.:
polymake: a framework for analyzing convex polytopes.
in {\em Polytopes --- Combinatorics and Computation}, edited by G Kalai and GM Ziegler (2000) 43--74.
\bibitem{Kimura1980}
Kimura, M.:
A simple method for estimating evolutionary rates of base substitution through comparative studies of nucleotide sequences.
Journal of Molecular Evolution {\bf 16} (1980) 111--120.
%
\bibitem {Neyman1971}
Neyman, J.:
Molecular studies of evolution: a source of novel statistical problems.
In {\em Statistical decision theory and related topics} edited by Gupta, S.,
Yackel, J., {New York Academic Press}, (1971) 1--27.
%
\bibitem {Jukes1969}
Jukes, H.T. and Cantor, C.:
Evolution of protein molecules.
In {\em Mammalian Protein Metabolism} edited by HN Munro, 
{New York Academic Press}, (1969) 21--32.
\bibitem{Levy2005}
Levy, D., Yoshida, R. and Pachter, L.:
Neighbor-joining with phylogenetic diversity estimates.
Molecular Biology and Evolution {\bf 23} (2006) 491--498.
\bibitem{Mihaescu2006}
Mihaescu, R., Levy, D., and Pachter, L.:
Why Neighbor-Joining Works. (2006) arXiv:cs.DS/0602041.
%
\bibitem{Olsen1994}
Olsen, G.J., Matsuda, H., Hagstrom, R., and Overbeek, R.:
fastDNAml: A tool for construction of phylogenetic trees of DNA sequences using maximum likelihood.
Comput. Appl. Biosci. {\bf 10} (1994) 41--48.

\bibitem{Ota2000}
Ota, S. and Li, WH.:
NJML: A Hybrid algorithm for the neighbor-joining and maximum likelihood methods.
Molecular Biology and Evolution {\bf 17} 9 (2000) 1401--1409.

\bibitem{Saitou1987}
Saitou, N. and Nei, M.:
The neighbor joining method: a new method for reconstructing phylogenetic trees. Molecular Biology and Evolution {\bf 4} (1987) 406--425.

\bibitem{Steel2006}
Gascuel, O. and Steel, M.: 
Neighbor-joining revealed. 
Molecular Biology and Evolution {\bf 23} (2006) 1997-­2000.

\bibitem{Studier1988}
Studier, J.A. and Keppler, K.J.:
A note on the neighbor-joining method of Saitou and Nei.
Molecular Biology and Evolution {\bf 5} (1988) 729--731.

\bibitem{Yang1997} Yang, Z: PAML: A program package for phylogenetic analysis by maximum likelihood. 
{\em CABIOS} {\bf 15} (1997) 555--556.

\bibitem {Yang2000}
Yang, Z.:
Complexity of the simplest phylogenetic estimation problem.
 Proceedings of the Royal Society B: Biological Sciences {\bf 267} (2000) 109-116.

\bibitem{Ziegler95}
Ziegler, G.:
{\em Lectures on Polytopes}. {Springer-Verlag} 1995.
%
\end{thebibliography}
